\documentclass{article}
\usepackage{amssymb,amsmath,theorem,euscript}

\input{epsf}

\newcounter{sec}

\def\sm{\smallskip}

\newcounter{punct}[sec]

\def\punct{\refstepcounter{punct}{\arabic{sec}.\arabic{punct}.  }}

\def\COUNTERS{\addtocounter{sec}{1}
              \setcounter{punct}{0}
          \setcounter{equation}{0}
          \setcounter{theorem}{0}
                  }

\newtheorem{theorem}{Theorem}[sec]
\newtheorem{proposition}[theorem]{Proposition}
\newtheorem{lemma}[theorem]{Lemma}
\newtheorem{corollary}[theorem]{Corollary}

\begin{document}

 \def\ov{\overline}
\def\wt{\widetilde}
 \newcommand{\rk}{\mathop {\mathrm {rk}}\nolimits}
\newcommand{\Aut}{\mathop {\mathrm {Aut}}\nolimits}
\newcommand{\Out}{\mathop {\mathrm {Out}}\nolimits}
 \newcommand{\tr}{\mathop {\mathrm {tr}}\nolimits}
  \newcommand{\diag}{\mathop {\mathrm {diag}}\nolimits}
  \newcommand{\supp}{\mathop {\mathrm {supp}}\nolimits}
  \newcommand{\indef}{\mathop {\mathrm {indef}}\nolimits}
  \newcommand{\dom}{\mathop {\mathrm {dom}}\nolimits}
  \newcommand{\im}{\mathop {\mathrm {im}}\nolimits}
 
\renewcommand{\Re}{\mathop {\mathrm {Re}}\nolimits}

\def\Br{\mathrm {Br}}

\def\O{\mathbb{O}}

\def\SL{\mathrm {SL}}
\def\SU{\mathrm {SU}}
\def\GL{\mathrm {GL}}
\def\U{\mathrm U}
\def\OO{\mathrm O}
 \def\Sp{\mathrm {Sp}}
 \def\SO{\mathrm {SO}}
\def\SOS{\mathrm {SO}^*}
 \def\Diff{\mathrm{Diff}}
 \def\Vect{\mathfrak{Vect}}
\def\PGL{\mathrm {PGL}}
\def\PU{\mathrm {PU}}
\def\PSL{\mathrm {PSL}}
\def\Symp{\mathrm{Symp}}
\def\End{\mathrm{End}}
\def\Mor{\mathrm{Mor}}
\def\Aut{\mathrm{Aut}}
 \def\PB{\mathrm{PB}}
 \def\cA{\mathcal A}
\def\cB{\mathcal B}
\def\cC{\mathcal C}
\def\cD{\mathcal D}
\def\cE{\mathcal E}
\def\cF{\mathcal F}
\def\cG{\mathcal G}
\def\cH{\mathcal H}
\def\cJ{\mathcal J}
\def\cI{\mathcal I}
\def\cK{\mathcal K}
 \def\cL{\mathcal L}
\def\cM{\mathcal M}
\def\cN{\mathcal N}
 \def\cO{\mathcal O}
\def\cP{\mathcal P}
\def\cQ{\mathcal Q}
\def\cR{\mathcal R}
\def\cS{\mathcal S}
\def\cT{\mathcal T}
\def\cU{\mathcal U}
\def\cV{\mathcal V}
 \def\cW{\mathcal W}
\def\cX{\mathcal X}
 \def\cY{\mathcal Y}
 \def\cZ{\mathcal Z}
\def\0{{\ov 0}}
 \def\1{{\ov 1}}
 \def\frA{\mathfrak A}
 \def\frB{\mathfrak B}
\def\frC{\mathfrak C}
\def\frD{\mathfrak D}
\def\frE{\mathfrak E}
\def\frF{\mathfrak F}
\def\frG{\mathfrak G}
\def\frH{\mathfrak H}
\def\frI{\mathfrak I}
 \def\frJ{\mathfrak J}
 \def\frK{\mathfrak K}
 \def\frL{\mathfrak L}
\def\frM{\mathfrak M}
 \def\frN{\mathfrak N} \def\frO{\mathfrak O} \def\frP{\mathfrak P} \def\frQ{\mathfrak Q} \def\frR{\mathfrak R}
 \def\frS{\mathfrak S} \def\frT{\mathfrak T} \def\frU{\mathfrak U} \def\frV{\mathfrak V} \def\frW{\mathfrak W}
 \def\frX{\mathfrak X} \def\frY{\mathfrak Y} \def\frZ{\mathfrak Z} \def\fra{\mathfrak a} \def\frb{\mathfrak b}
 \def\frc{\mathfrak c} \def\frd{\mathfrak d} \def\fre{\mathfrak e} \def\frf{\mathfrak f} \def\frg{\mathfrak g}
 \def\frh{\mathfrak h} \def\fri{\mathfrak i} \def\frj{\mathfrak j} \def\frk{\mathfrak k} \def\frl{\mathfrak l}
 \def\frm{\mathfrak m} \def\frn{\mathfrak n} \def\fro{\mathfrak o} \def\frp{\mathfrak p} \def\frq{\mathfrak q}
 \def\frr{\mathfrak r} \def\frs{\mathfrak s} \def\frt{\mathfrak t} \def\fru{\mathfrak u} \def\frv{\mathfrak v}
 \def\frw{\mathfrak w} \def\frx{\mathfrak x} \def\fry{\mathfrak y} \def\frz{\mathfrak z} \def\frsp{\mathfrak{sp}}
 \def\bfa{\mathbf a} \def\bfb{\mathbf b} \def\bfc{\mathbf c} \def\bfd{\mathbf d} \def\bfe{\mathbf e} \def\bff{\mathbf f}
 \def\bfg{\mathbf g} \def\bfh{\mathbf h} \def\bfi{\mathbf i} \def\bfj{\mathbf j} \def\bfk{\mathbf k} \def\bfl{\mathbf l}
 \def\bfm{\mathbf m} \def\bfn{\mathbf n} \def\bfo{\mathbf o} \def\bfp{\mathbf p} \def\bfq{\mathbf q} \def\bfr{\mathbf r}
 \def\bfs{\mathbf s} \def\bft{\mathbf t} \def\bfu{\mathbf u} \def\bfv{\mathbf v} \def\bfw{\mathbf w} \def\bfx{\mathbf x}
 \def\bfy{\mathbf y} \def\bfz{\mathbf z} \def\bfA{\mathbf A} \def\bfB{\mathbf B} \def\bfC{\mathbf C} \def\bfD{\mathbf D}
 \def\bfE{\mathbf E} \def\bfF{\mathbf F} \def\bfG{\mathbf G} \def\bfH{\mathbf H} \def\bfI{\mathbf I} \def\bfJ{\mathbf J}
 \def\bfK{\mathbf K} \def\bfL{\mathbf L} \def\bfM{\mathbf M} \def\bfN{\mathbf N} \def\bfO{\mathbf O} \def\bfP{\mathbf P}
 \def\bfQ{\mathbf Q} \def\bfR{\mathbf R} \def\bfS{\mathbf S} \def\bfT{\mathbf T} \def\bfU{\mathbf U} \def\bfV{\mathbf V}
 \def\bfW{\mathbf W} \def\bfX{\mathbf X} \def\bfY{\mathbf Y} \def\bfZ{\mathbf Z} \def\bfw{\mathbf w}
 \def\R {{\mathbb R }} \def\C {{\mathbb C }} \def\Z{{\mathbb Z}} \def\H{{\mathbb H}} \def\K{{\mathbb K}}
 \def\N{{\mathbb N}} \def\Q{{\mathbb Q}} \def\A{{\mathbb A}} \def\T{\mathbb T} \def\P{\mathbb P} \def\G{\mathbb G}
 \def\bbA{\mathbb A} \def\bbB{\mathbb B} \def\bbD{\mathbb D} \def\bbE{\mathbb E} \def\bbF{\mathbb F} \def\bbG{\mathbb G}
 \def\bbI{\mathbb I} \def\bbJ{\mathbb J} \def\bbK{\mathbb K} \def\bbL{\mathbb L} \def\bbM{\mathbb M} \def\bbN{\mathbb N} \def\bbO{\mathbb O}
 \def\bbP{\mathbb P} \def\bbQ{\mathbb Q} \def\bbS{\mathbb S} \def\bbT{\mathbb T} \def\bbU{\mathbb U} \def\bbV{\mathbb V}
 \def\bbW{\mathbb W} \def\bbX{\mathbb X} \def\bbY{\mathbb Y} \def\kappa{\varkappa} \def\epsilon{\varepsilon}
 \def\phi{\varphi} \def\le{\leqslant} \def\ge{\geqslant}

\def\UU{\bbU}
\def\Mat{\mathrm{Mat}}
\def\tto{\rightrightarrows}

\def\Gr{\mathrm{Gr}}

\def\graph{\mathrm{graph}}

\def\la{\langle}
\def\ra{\rangle}

\def\B{\mathrm B}
\def\Int{\mathrm{Int}}
\def\LGr{\mathrm{LGr}}


\def\I{\mathbb I}
\def\M{\mathbb M}
\def\T{\mathbb T}
\def\S{\mathrm S}

\def\Lat{\mathrm{Lat}}
\def\BT{\mathrm{BT}}

\def\LLat{\mathrm{LLat}} 
\def\Mod{\mathrm{Mod}}
\def\LMod{\mathrm{LMod}}
\def\Naz{\mathrm{Naz}}
\def\naz{\mathrm{naz}}
\def\bNaz{\mathbf{Naz}}
\def\AMod{\mathrm{AMod}}
\def\ALat{\mathrm{ALat}}
\def\Coll{\mathrm{Coll}}

\def\Ver{\mathrm{Vert}}
\def\Bd{\mathrm{Bd}}
\def\We{\mathrm{We}}
\def\Heis{\mathrm{Heis}}

\def\bbot{{\bot\!\!\!\bot}}
\def\frcoll{\mathfrak{Coll}}

\begin{center}
\Large\bf

On $p$-adic colligations and 'rational maps' of  Bruhat-Tits trees

\bigskip

\large\sc
Yury A. Neretin\footnote{Supported by the grants FWF, P22122, P25142.}
\end{center}

{\small Consider   matrices of order $k+N$  over $p$-adic field determined up to
conjugations by elements of $GL$ over  $p$-adic integers. We define a product of such conjugacy classes
and construct the analog of characteristic functions (transfer functions), they are maps
from Bruhat-Tits trees to Bruhat-Tits buildings. We also examine categorical quotient for usual 
operator colligations.}

\section{Introduction}

\COUNTERS

{\bf\punct Notation.} Denote by $1=1_\alpha$ the unit matrix of order $\alpha$.
 Below $K$ is an  infinite  field%
 \footnote{We prefer infinite fields, otherwise the rational function 
 (\ref{eq:char-raz}) is not well defined.}, $\K$ is a locally compact non-Archimedian
 field, $\O  \subset \K$ is the ring of integers. In both cases we keep in mind the $p$-adic fields.
Let $\Mat(n)=\Mat(n,K)$ be
the space of matrices of order $n$ over $K$, $\GL(n,K)$ the group of invertible matrices
of order $n$. We say that an $\infty\times\infty$ matrix $g$ is {\it finite}
if $g-1$ has finite number of nonzero matrix elements%
\footnote{Thus $1_\infty$ is finite and 0 is not finite}.
Denote by $\Mat(\infty)=\Mat(\infty,\K)$  the space of finite 
$\infty\times\infty$ matrices,
by $\GL(\infty,K)$ the group of finite invertible finite matrices.

\sm

{\bf\punct Colligations.} Consider the space $\Mat(\alpha+\infty,K)$
of finite block complex matrices
$g=\begin{pmatrix}a&b\\c&d\end{pmatrix}$
 of size $(\alpha+\infty)\times(\alpha+\infty)$.
 Represent the group $\GL(\infty,K)$ as the group of matrices
 of the form $\begin{pmatrix} 1_\alpha&0\\0&u\end{pmatrix}$ of size $\alpha+\infty$.
 Consider conjugacy classes of
 $\Mat(\alpha+\infty,K)$ with respect to $\GL(\infty,K)$, i.e.,
 matrices determined up to the equivalence 
 \begin{equation}
 \begin{pmatrix}a&b\\c&d \end{pmatrix}
 \sim
 \begin{pmatrix} 1_\alpha&0\\0&u\end{pmatrix}
 \begin{pmatrix}a&b\\c&d\end{pmatrix}
 \begin{pmatrix} 1_\alpha&0\\0&u\end{pmatrix}^{-1},\qquad \text{where $u\in\GL(\infty,\K)$}
 \label{eq:sim}
 .\end{equation}
 We call conjugacy classes by {\it colligations} (another term is '{\it nodes}').
 Denote by 
 $$\Coll(\alpha)=\Coll(\alpha,K)$$
  the set of equivalence classes.
There is a natural multiplication on $\Coll(\alpha,K)$,
it is given by
\begin{equation}
 \begin{pmatrix}a&b\\c&d\end{pmatrix}\circ
  \begin{pmatrix}p&q\\r&t\end{pmatrix}
  =
   \begin{pmatrix}a&b&0\\c&d&0\\ 0&0&1\end{pmatrix}
    \begin{pmatrix}p&0&q\\0&1&0\\ r&0&t\end{pmatrix} 
    =
    \begin{pmatrix}
    ap&b&aq\\cp&d&cq\\r&0&t   
    \end{pmatrix}
    \label{eq:circ}
.\end{equation}
The size of the last matrix is
$$
\alpha+\infty+\infty=\alpha+\infty
.$$

The following statement is straightforward.

\begin{proposition}
\label{pr:circ-1}
{\rm a)} The $\circ$-multiplication is a well-defined operation
$$
\Coll(\alpha)\times\Coll(\alpha)\to\Coll(\alpha)
.
$$

\sm

{\rm b)} The $\circ$-multiplication is associative.
\end{proposition} 

There is a way to visualize this multiplication. We write the 
following 'perverse' equation for eigenvalues:
\begin{equation}
\begin{pmatrix}q\\ x \end{pmatrix}
=
 \begin{pmatrix}a&b\\c&d\end{pmatrix}
 \begin{pmatrix}p\\ \lambda x \end{pmatrix}
 ,
 \label{eq:perverse}
\end{equation}
where $\lambda\in K$. Equivalently,
\begin{align}
q=ap+\lambda bx;
\label{eq:elimination-1}
\\
x=cp+\lambda dx.
\label{eq:elimination-2}
\end{align}
We express $x$ from (\ref{eq:elimination-2}), 
$$
x=(1-\lambda d)^{-1} cp,
$$
substitute it to (\ref{eq:elimination-1}), and get 
$$
q=\chi_g(\lambda) p
,$$
where $\chi_g(\lambda)$ 
\begin{equation}
\chi_g(\lambda)=a+\lambda b(1-\lambda d)^{-1}c
\label{eq:char-raz}
\end{equation}
is a rational function $K\to \Mat(\alpha)$.
It is called {\it characteristic function of $g$'}.
The following statement is obvious.

\begin{proposition}
\label{pr:independence}
If $g_1$ and $g_2$ are contained in the same conjugacy class,
then their characteristic functions coincide.
\end{proposition}

The next statement can be verified by a straightforward calculation (for a more reasonable
proof, see below Theorem \ref{th:product}).

\begin{theorem}
\label{th:circ-char}
$$
\chi_{g\circ h}(\lambda)=\chi_g(\lambda)\chi_h(\lambda)
.
$$
\end{theorem}

\begin{theorem}
\label{th:predst}
Let $K$ be algebraically closed. 
Then any rational map $K\to\Mat(\alpha,K)$ regular at $0$ has the form 
$\chi_g(\lambda)$ for a certain $g\in\Mat(\infty, K)$.
\end{theorem}

See. e.g., \cite{Dym}, Theorem 19.1.

\sm


{\bf\punct Origins of the colligations.} The colligations and
the characteristic functions
appeared independently in  spectral theory of non-self-adjoint 
operators (M.S.Livshits, 1946) and in system theory, see e.g., \cite{Liv1}, \cite{Liv2}, \cite{Pot},
\cite{MH}, \cite{Haz}, \cite{RR}, 
 \cite{Dym}, \cite{Bro}, \cite{GGK}. It seems that in both cases there are no visible reason
 to pass to $p$-adic case.

However, colligations and colligation-like objects arose by independent reasons
 in representation theory of
infinite-dimensional classical groups, see \cite{Olsh-GB}, \cite{Ner-book}.

First, consider a locally compact non-Archimedian field $\K$ and
the double cosets 
$$M=\SL(2,\O)\setminus \SL(2,\K)/\SL(2,\O).$$
 The space of functions on $M$
is a commutative algebra with respect to the convolution on $\SL(2 ,\K)$.
This algebra acts in the space of $\SL(2,\O)$-fixed vectors of any unitary representation of
$\SL(2,\K)$.  Next (see Ismagilov \cite{Ism1}, 1967),
let us replace $\K$ by a non-Archimedian non-locally compact field 
(i.e., the residue field is infinite \cite{Ism1} or
the  norm group is non-discrete \cite{Ism3}). Then there is no convolution,
however double cosets have a natural structure of a semigroup, and this semigroup acts in the space
of $\SL(2,\O)$-fixed vectors of any unitary representation of $\SL(2,\K)$. In particular, this allows to classify all irreducible unitary representations of $\SL(2,\K)$ having a non-zero
$\SL(2,\O)$-fixed vector.

It appeared that these phenomena (semigroup structure on double cosets $L\setminus G/L$
 for infinite dimension groups%
 \footnote{There is a elementary explanation initially proposed by Olshanski: such semigroups are limits
 of Hecke-type algebras at infinity. For more details,
 see \cite{Ner-haar}}%
 $^,$%
 \footnote{Conjugacy classes are special cases of double cosets, conjugacy classes
 $G$ with respect to $L$ are double cosets $L\setminus (G\times L)/L$, where $L$ is embedded to
 $G\times L$ diagonally, $l\mapsto  (l,l)$.}
  and
actions of this semigroup in the space of $L$-fixed vectors)
are quite general, see, e.g., \cite{Olsh-tree}, \cite{Olsh-GB}, \cite{Olsh-topics}, \cite{Ner-book}, \cite{Ner-faa},
\cite{Ner-char}.

 In \cite{Ner-buildings} there was proposed a way to construct representations of infinite-dimensional 
 $p$-adic groups, in particular there appeared  semigroups of double cosets and $p$-adic colligation-like structures. The present work is a simplified parallel of  \cite{Ner-buildings}. If we  look to the equivalence
 (\ref{eq:sim}), then a $p$-adic field is an representative of non-algebraically closed fields. However,
 \cite{Ner-buildings} suggests another equivalence,
 \begin{equation}
 \begin{pmatrix}a&b\\c&d \end{pmatrix}
 \sim
 \begin{pmatrix} 1_\alpha&0\\0&u\end{pmatrix}
 \begin{pmatrix}a&b\\c&d\end{pmatrix}
 \begin{pmatrix} 1_\alpha&0\\0&u\end{pmatrix}^{-1},\qquad \text{where $u\in\GL(\infty,\O)$}
 \label{eq:sim2}
 \end{equation} 
 (we conjugate by the group $\GL(\infty,\O)$ of integer matrices).
 Below we construct analogs of characteristic functions for this equivalence and  get 'rational' maps
 from Bruhat-Tits trees to Bruhat-Tits buildings (for $\alpha=1$ we get  maps from trees to trees),
 the characteristic function (\ref{eq:char-raz}) is its boundary values on the absolute of the tree.
 
 It is interesting that maps of this type arise in theory of Berkovich rigid analytic spaces%
 \footnote{In Berkovich theory objects are lager than trees and buildings. However, our 'characteristic functions' admit extensions to these larger objects},
 see \cite{Bak}, \cite{BR}, \cite{Con}. However I do not understand links between two  points of view.
 For instance, we show that any rational map of a projective line $\P\Q_p^1$ to itself 
 admits a continuation to the Bruhat-Tits tree, and such continuations are enumerated
 by the set $\GL(\infty,\Q_p)/\GL(\infty,\O)$. In Berkovich theory continuations of this type are canonical.


\sm

{\bf \punct Structure of the paper.} In Section 2 we consider characteristic functions over algebraically closed field. 
We discuss categorical quotient $[\Coll(\alpha)]$ of $\Coll(\alpha)$ with respect to the equivalence (\ref{eq:sim}),
the main statement here is Theorem \ref{th:isomorphism}. 
In Section 3 we examine the case $\alpha=1$. We show that the semigroup $[\Coll(1)]$ is commutative.
Also we show that for non-algebraically closed field any rational function 
$K\to K$ is a characteristic function.

In Section 4 we consider $p$-adic fields and introduce characteristic functions for conjugacy classes
of $\GL(\alpha+\infty,\Q_p)$ by $\GL(\infty,\O_p)$. 

In Section 5 we briefly discuss conjugacy classes of $\GL(\alpha+m\infty,\Q_p)$ with respect to
$\GL(\infty,\O_p)$.


\section{Formalities. Algebraically closed fields}

\COUNTERS

In this section $K$ is an algebraically closed field. For exposition 
of basic classical theory, see the textbook of Dym \cite{Dym}, Chapter 19. See more
in \cite{MH}, \cite{Haz}, \cite{RR}, \cite{Sot}. Our 'new' element
is the categorical quotient%
\footnote{I have not met discussion of this topic,
however sets of 'nonsingular points' of $\Coll(\alpha)$ and its completions were discussed in literature,
see \cite{Haz}, \cite{RR}.}
 (in a wider generality it was discussed 
in \cite{Ner-invariant}). 

\sm

Denote by $\P K^1$ the projective line over $K$. For an even dimensional linear
space $W$ denote by $\Gr(W)$ the Grassmannian of subspaces of dimension $\frac 12 \dim V$.

\sm


{\bf\punct Colligations.}
Fix $\alpha\ge 0$, $N>0$. Consider the space of matrices $\Mat(\alpha+N,K)$,
we write its elements as block matrices $g=\begin{pmatrix}a&b\\c&d\end{pmatrix}$.
Consider the group $\GL(N,K)$, we represent its elements as block matrices
$\begin{pmatrix} 1_\alpha&0\\0&u\end{pmatrix}$. Denote by
$\Coll_N(\alpha,K)$ the space of conjugacy classes of
$\Mat(\alpha+N)$ with respect to $\GL(N,K)$, see (\ref{eq:sim}).
Denote by $[\Coll_N(\alpha,K)]$ the categorical quotient (see, e.g., \cite{PV}), i.e.,
the spectrum of the algebra of $\GL(N,K)$-invariant 
polynomials on $\Mat(\alpha+N)$.

 
\sm

{\bf\punct Characteristic function.}
For an element $g=\begin{pmatrix}a&b\\c&d\end{pmatrix}$ of $\Mat(\alpha+N)$
we assign the characteristic function
\begin{equation}
\chi_g(\lambda) = a+\lambda b(1-\lambda d)^{-1} c,\qquad \text{$\lambda$ ranges in $K$.}
\label{eq:char-2}
\end{equation}

If $d$ is invertible, we extend this function to the point $\lambda=\infty$
by setting
$$
\chi_g(\infty)=a-bd^{-1}c.
$$
Passing to the coordinate $s=\lambda^{-1}$ on $\P K^1$, we get
$$
\chi_g(s)=a+b(s-d)^{-1} c.
$$

\begin{theorem}
Any rational function $K\to\Mat(\alpha,K)$ regular at $0$
is a characteristic function of an operator colligation. 
\end{theorem}

See, e.g., \cite{Dym}, Theorem 19.1.


\sm

{\bf \punct The characteristic function as a map $\P K^1\to \Gr(K^{2\alpha})$.}
See \cite{MH}, \cite{Haz}, \cite{Sot}.
If $\lambda_0$ is a regular point of $\chi_g(\lambda)$, we consider its graph $\cX_g(\lambda_0)$,
$$
\cX_g(\lambda_0) \subset K^\alpha\oplus K^\alpha
.
$$
Singularities of rational maps of $\P K^1$ to a projective variety
$\Gr(K^\alpha\oplus K^\alpha)$ are removable. Let us remove a singularity explicitly
at a pole $\lambda=\lambda_0$.
We can represent $\chi_g(\lambda)$ as
\begin{equation}
A(\lambda-\lambda_0)
\begin{pmatrix}\frac{h_1}{(\lambda-\lambda_0)^{m_1}}&0&\dots\\
0&\frac{h_2}{(\lambda-\lambda_0)^{m_2}}&\dots\\
\vdots& \vdots & \ddots
 \end{pmatrix} B(\lambda-\lambda_0)+ S(\lambda-\lambda_0)
,
\end{equation}
where $A(\dots)$, $B(\dots)$ are polynomial functions $K\to\Mat(\alpha)$,
$A(0)$, $B(0)$ are invertible,
the exponents $m_i$ satisfy $m_1\ge m_2\ge \dots$, and $S(\lambda)$ is a rational functions
$K\to\Mat(\alpha)$ having zero of any prescribed  order $M>0$ (proof of this
is a straightforward repetition   of the Gauss elimination procedure).
Denote by $e_j$ the standard basis in $K^\alpha$. Consider the 
subspace $L$ in $K^\alpha\oplus K^\alpha$ generated by vectors
\begin{align*}
e_i\oplus 0,\qquad \text{for $m_i>0$};
\\
h_j e_j \oplus e_j,\qquad \text{for $m_j=0$};
\\
0\oplus e_l,\qquad \text{for $m_l<0$}.
\end{align*}
Applying the operator $A(0)\oplus B(0)^{-1}$ to $L$ we get $\chi_g(\lambda_0)$.

\sm

{\bf \punct An exceptional divisor.} A characteristic function is not sufficient for a reconstruction
of a colligation.
Indeed, 
consider a block matrix of size $\alpha+k+l$
$$
\begin{pmatrix}
a&b&0\\c&d&0\\0&0&e
\end{pmatrix}
$$
Then characteristic the function is independent on $e$.

 For $g\in\Coll_N(\alpha)$ we define an additional invariant, a divisor%
\footnote{i.e. a finite  set with multiplicities.}
 $\Xi_g\subset \P K^1$ in the following way: 
 $\Xi_g$  as the divisor  of zeros 
 of the polynomial
 $$
 p_g(\lambda)=\det(1-\lambda d)
 $$
plus $\lambda=\infty$ with multiplicity $N-\deg p_g$.
In the coordinate $s=\lambda^{-1}$ this divisor is simply the set
of eigenvalues of $d$. 

\begin{proposition}
\label{pr:det}
$$
\det\chi_g(\lambda)=\frac{\det\begin{pmatrix}a&-\lambda b\\ c&1-\lambda d\end{pmatrix}}
                        {\det(1-\lambda d)}
                        .
                        $$
\end{proposition}

{\sc Proof.} We apply the formula for the determinant of a block matrix.
\hfill $\square$

\sm

\begin{corollary}
The divisor $\Xi_g$ contains the divisor of poles of $\det\chi_g(\lambda)$.
\end{corollary}

\begin{theorem}
For any rational function $K\to \Mat(\alpha)$ regular at $0$ there is a  colligation $g$
such that the divisor  $\Xi_g$ coincides with the divisors of poles of $\det\chi_g(\lambda)$.
\end{theorem}

See \cite{Dym}, Theorem 19.8. Such colligations $g$ are called {\it minimal}.

\sm


{\bf \punct Invariants.}

\begin{theorem}
A point $\frg$ of the categorical quotient $[\Coll_N(\alpha)]$
is uniquely determined by the characteristic function $\chi_\frg(\lambda)$ and the
divisor $\Xi_\frg$.
\end{theorem}

{\sc Proof.} Let us describe $\GL(N,K)$-invariants on $\Mat(\alpha+N)$.
A point $g=\begin{pmatrix}a&b\\c&d \end{pmatrix}\in \Mat(\alpha+N)$
can be regarded as the following collection of data:

a) the matrix $d$;

b) $\alpha$ vectors (columns $c[j]$ of $c$);

c) $\alpha$ covectors (rows $b[i]$ of $b$);

d) scalars $a_{ij}$.

\sm

The algebra of invariants (see \cite{Pro}, Section 11.8.1) is generated by the following polynomials
\begin{align}
&b[i] d^kc[j],
\label{eq:invariant-1}
\\
&\tr d^k,
\label{eq:invariant-2}
\\
&a_{ij}.
\label{eq:invariant-3}
\end{align}
Expanding the characteristic function in $\lambda$, 
$$
\chi_g(\lambda)=a+\sum_{k=0}^\infty \lambda^{k+1} bd^k c
$$
we get in coefficients all the invariants (\ref{eq:invariant-1}),
(\ref{eq:invariant-3}). Expanding
$$
\ln p_g(\lambda)=\ln \det(1-\lambda d)=-\sum_{j=k}^\infty \frac 1k \lambda^k \tr d^k
,
$$  
we get all invariants (\ref{eq:invariant-2}).
\hfill $\square$

\begin{corollary}
Any point of $[\Coll_N(\alpha)]$ has a  representative of the form
$$
\begin{pmatrix}
a&b&0\\
c&d&0\\
0&0&e
\end{pmatrix}
,$$
where $e$ is diagonal matrix and the colligation $\begin{pmatrix}a&b\\c&d\end{pmatrix}$ is minimal.
\end{corollary}


{\bf\punct $\circ$-product.} Now we define the operation 
$$
\Coll_{N_1}(\alpha)\times \Coll_{N_2}(\alpha)\to \Coll_{N_1+N_2}(\alpha)
$$
 by the formula (\ref{eq:circ}).
 
 \begin{theorem}
 \label{th:prod-1}
 $$
 \mathrm {a)} \quad\quad
\chi_{g\circ h}(\lambda)=\chi_g(\lambda) \chi_h(\lambda).
$$
$$ 
\mathrm {b)} \quad\quad \Xi_{g\circ h}=\Xi_g+\Xi_h.
$$
 \end{theorem}
 
  The statement b) is obvious, a) is well-known (see a proof below, Theorem \ref{th:product}).
  
 \begin{corollary}
 The $\circ$-multiplication is well defined as an operation on categorical
 quotients,
$$
[\Coll_{N_1}(\alpha)]\times [\Coll_{N_2}(\alpha)]\to [\Coll_{N_1+N_2}(\alpha)]
$$ 
 \end{corollary}

{\sc Proof.} Indeed, invariants of $g\circ h$ are determined by invariants 
of $g$ and $h$.\hfill $\square$

\sm


{\bf\punct The space $\Coll_\infty(\alpha)$.}
Consider the natural map 
$$I_N:=\Mat(\alpha+N)\to \Mat(\alpha+N+1)$$
defined by
$$
I_N:
\begin{pmatrix} a&b\\c&d\end{pmatrix}
\mapsto 
\begin{pmatrix} a&b&0\\c&d&0\\ 0&0&1\end{pmatrix}
.
$$
We have 
\begin{align*}
\chi_{I_N g}(\lambda)=\chi_g(\lambda);
\\
\Xi_{I_Ng}=\Xi_g+\{1\}
,
\end{align*}
where $\{1\}$ is the point $1\in K$.

\begin{lemma}
The induced map $[\Coll_N(\alpha)]\to [\Coll_{N+1}(\alpha)]$
is an embedding.
\end{lemma}

{\sc Proof.} The restriction of invariants (\ref{eq:invariant-1})--(\ref{eq:invariant-3})
defined on $\Mat(\alpha+N+1)$ to the subspace
 $\Mat(\alpha+N)$ gives the same expressions for $\Mat(\alpha+N)$.
\hfill $\square$

\sm

Thus, we can define a space $[\Coll(\alpha)]=[\Coll_\infty(\alpha)]$ as an inductive limit 
$$
[\Coll_\infty(\alpha)]=
\lim_{N\to\infty}[\Coll_N(\alpha)]
.$$ 
It is equipped with the associative $\circ$-multiplication.

Characteristic function of an element $g\in[\Coll_\infty(\alpha)]$ can be defined in two equivalent ways.
The first way, we  write the expression (\ref{eq:char-2}) for infinite matrix $g$.
The second way. We choose large $N$ such that $g$ has the following representation
\begin{equation}
g=
\left(
\begin{array}{cccl}
a&b&0&\}\alpha\\
c&d&0&\}N\\
0&0&1_\infty&\}\infty
\end{array}
\!\!\!\!\!\!\!\!\!\!\!\!\!\!\!\!\!\!
\right)
\label{eq:cut}
\end{equation}
and write the characteristic function for the upper left block of the size $\alpha+N$.

Next, we define the {\it exceptional divisor} $\Xi_g$ in $\P K^1$.
 We represent $g$ in the form (\ref{eq:cut}),
write the exceptional divisor for $\begin{pmatrix}a&b\\c&d \end{pmatrix}$,
and add the point $\lambda=1$ with multiplicity $\infty$ (in particular, the multiplicity
of $1$ always is infinity).
Thus, we can regard the 'divisor' as  a function 
$$
\xi:\P K^1\to \Z_+\cup \infty
$$
satisfying the following condition.

\sm

a) $\xi(\lambda)=0$ for all but finite number of $\lambda$.

\sm

b) $\xi(1)=\infty$, at all other points $\xi$ is finite.

\sm

We reformulate  statements obtained above in the following form.
Denote by $\Gamma_\alpha$ the semigroup  of rational maps $\chi:\P K^1\to \Mat(\alpha)$
regular at the point $\lambda=0$. Denote by $\Delta$ the set 
of all divisors in the sense described above.
We equip $\Delta$ with the operation of addition.

Next, consider the subsemigroup  
$R_\alpha\subset\Gamma\times \Delta$ consisting of pairs
$(\chi,\Xi)$ such that divisor of the denominator
of $\det\chi(\lambda)$ is contained in the divisor $\Xi$.

\begin{theorem}
\label{th:isomorphism}
 The map $g\mapsto (\chi_g,\Xi_g)$
is an isomorphism of semigroups $[\Coll_\infty(\alpha)]$
and $R_\alpha$.
\end{theorem}

Notice, that the semigroup $[\Coll_\infty(\alpha)$ itself is not a product 
of semigroup of characteristic functions and an Abelian semigroup.
A similar object appeared in \cite{Ner-book}, IX.2.


\section{The case $\alpha=1$}

\COUNTERS

{\bf \punct Commutativity.}

\begin{theorem}
The semigroup $[\Coll_\infty(1)]$ is commutative.
\end{theorem}

{\sc Proof.} Indeed, $\Gamma_1$ is commutative,  therefore 
$\Gamma_1\times \Delta$  is commutative.
\hfill $\square$

\sm

{\sc Remark.} The semigroup $\Coll_\infty(1)$ is not commutative,
\begin{align*}
\begin{pmatrix}
1&0\\1&1
\end{pmatrix}
\circ
\begin{pmatrix}
1&1\\0&1
\end{pmatrix}
=\begin{pmatrix}
1&0&0\\
1&1&1\\
0&0&1
\end{pmatrix},
\\
\begin{pmatrix}
1&1\\0&1
\end{pmatrix}
\circ
\begin{pmatrix}
1&0\\1&1
\end{pmatrix}
=\begin{pmatrix}
1&1&0\\
0&1&0\\
1&0&1
\end{pmatrix}
,\end{align*}
and blocks '$d$' in the right hand side have different Jordan forms.
\hfill $\square$

\sm

{\bf \punct Commutativity. Straightforward proof.}
The proof given below is not necessary in the contexts of this paper.
However, it shows that the commutativity  in certain sense is a non-obvious fact
(in particular, this proof can be modified for proofs of non-commutativity 
of $\circ$-products in some cases discussed in \cite{Ner-faa}).

First, an element of $\Coll_N(1)$ in a general position
can be reduced by a conjugation to the form
$$
\begin{pmatrix}
a&b_1&b_2&b_3&\dots
\\
c_1&\lambda_1&0&0&\dots\\
c_2&0&\lambda_2&0&\dots\\
c_3&0&0&\lambda_3&\dots\\
\vdots&\vdots&\vdots&\vdots &\ddots
\end{pmatrix}
,$$
where $\lambda_j$ are pairwise distinct.
To be short set $N=2$. Consider two matrices
$$
g=
\begin{pmatrix}
p&b_1&b_2
\\
c_1&\lambda_1&0\\
c_2&0&\lambda_2
\end{pmatrix}
\qquad
h=
\begin{pmatrix}
p&q_1&q_2
\\
r_1&\mu_1&0\\
r_2&0&\mu_2
\end{pmatrix}
$$
with $\lambda_1$, $\lambda_2$, $\mu_1$, $\mu_2$ being pairwise distinct. We
evaluate
$$
S=g\circ h =
\begin{pmatrix}
ap&b_1&b_2&aq_1&aq_2\\
c_1p&\lambda_1&0&c_1q_1&c_1q_2\\
c_2p&0&\lambda_2&c_2q_1&c_2q_2\\
r_1&0&0&\mu_1&0\\
r_2&0&0&0&\mu_2&
\end{pmatrix}
$$
and
\begin{multline*}
T= h\circ g=
\begin{pmatrix}
p&0& 0& q_1&q_2\\
0&1&  0&0&0\\
0&0& 1& 0&0\\
r_1&0&0& \mu_1&0\\
r_2&0& 0&0&\mu_2
\end{pmatrix}
\begin{pmatrix}
a&b_1&   b_2&       0&0\\
c_1&\lambda_1&  0&0&0\\
c_2&0&          \lambda_2&0&0\\
0&0&       0&     1&0\\
0&0&            0&0&1
\end{pmatrix}
=\\=
\begin{pmatrix}
ap&b_1p&b_2p&q_1&q_2\\
c_1&\lambda_1&0&0&0\\
c_2&0&\lambda_2&0&0\\
ar_1&b_1r_1&b_2r_1&\mu_1&0\\
ar_2&b_1r_2&b_2r_2&0&\mu_2
\end{pmatrix}
\end{multline*}

\begin{proposition}
In this notation,
$$
T=\begin{pmatrix}
1&0\\0&U
\end{pmatrix}^{-1} S
\begin{pmatrix}
1&0\\0&U
\end{pmatrix}
,$$
where 
$$
U=U_+^{-1} U_d U_-,
$$
matrices $U_+$, $U_-$ are upper {\rm(}lower{\rm)} triangular respectively,
$$
U_+= 
\begin{pmatrix}
1&0&\frac{c_1q_1}{\lambda_1-\mu_1}&\frac{c_1q_2}{\lambda_1-\mu_2}\\
0&1&\frac{c_2q_1}{\lambda_2-\mu_1}&\frac{c_2q_2}{\lambda_2-\mu_2}\\
0&0&1&0\\
0&0&0&1
\end{pmatrix},
\qquad
U_-=
\begin{pmatrix}
1&0&0&0\\
0&1&0&0\\
\frac{b_1r_1}{\lambda_1-\mu_1}&\frac{b_2r_1}{\lambda_2-\mu_1}&1&0\\
\frac{b_1r_2}{\lambda_1-\mu_2}&\frac{b_2r_2}{\lambda_2-\mu_2}&0&1
\end{pmatrix}
,
$$
and $U_d$ is a diagonal matrix with entries
\begin{align*}
p+\frac{q_1 r_1}{\lambda_1-\mu_1}+ \frac{q_2 r_2}{\lambda_1-\mu_2},\quad
p+\frac{q_1 r_1}{\lambda_2-\mu_1}+ \frac{q_2 r_2}{\lambda_2-\mu_2},
\\
\left(a+\frac{b_1c_1}{\mu_1-\lambda_1}+\frac{b_2c_2}{\mu_1-\lambda_2}\right)^{-1},\quad
\left(a+\frac{b_1c_1}{\mu_2-\lambda_1}+\frac{b_2c_2}{\mu_2-\lambda_2}\right)^{-1}
\end{align*}
\end{proposition}

{\sc Proof.} We represent $T$ and $S$ as block matrices,
$$
T=\begin{pmatrix}T_{11}&T_{12}\\T_{21}&T_{22}  \end{pmatrix},
\qquad
S=\begin{pmatrix}S_{11}&S_{12}\\S_{21}&S_{22}  \end{pmatrix}
$$
of size $(1+4)\times (1+4)$. We must verify equalities
\begin{align}
&UT_{22}=S_{22}U,
\label{eq:22}
\\
& UT_{21}=S_{21},\qquad T_{12}=S_{12}U
 .
 \label{eq:21}
\end{align}
Represent the first equality in the form 
\begin{equation}
U_d(U_-T_{22} U_-^{-1})=(U_+S_{22} U_+^{-1}) U_d
\label{eq:diag}
.\end{equation}
The matrices $U_\pm$ are chosen in such a way that
$$U_-T_{22} U_-^{-1}=U_+S_{22} U_+^{-1}
=\begin{pmatrix}\lambda_1&0&0&0\\
0&\lambda_2&0&0\\
0&0&\mu_1&0\\
0&0&0&\mu_2
 \end{pmatrix} 
$$
Therefore (\ref{eq:diag}) holds for any diagonal matrix $U_d$.
It remains to choose $U_d$ to satisfy (\ref{eq:21}).
\hfill $\square$

\sm

{\bf\punct  Linear-fractional transformations.}

\begin{proposition}
\label{pr:lin-frac}
Let $\begin{pmatrix}\alpha&\beta\\\gamma&\delta \end{pmatrix}$
be a nondegenerate $2\times 2$ matrix. Let $g\in \Coll_N(\alpha)$.
Let $\gamma g+\delta$ be nondegenerate%
\footnote{this is independent on the choice of a representative.}.
Let $\chi_g(s)$ be the characteristic function written in the coordinate
$s=\lambda^{-1}$ 
Then  the characteristic function of the colligation
$$
h=(\alpha g+\beta)(\gamma g+\delta)^{-1}
$$
is
$$
\left(\alpha \chi_g\left(\frac{\alpha s+\beta}{\gamma s+\delta}\right)+\beta\right)
\left(\gamma \chi_g\left(\frac{\alpha s+\beta}{\gamma s+\delta}\right)+\delta\right)^{-1}
.
$$
\end{proposition}

{\sc Proof.} We represent the equation
$$
\begin{pmatrix}q\\sx \end{pmatrix}=h \begin{pmatrix}p\\x \end{pmatrix}
$$
as 
$$
\begin{pmatrix} 
\alpha q+\beta p\\
\alpha sx+\beta x
\end{pmatrix}
=
\begin{pmatrix}
a&b\\c&d
\end{pmatrix}
\begin{pmatrix}
\gamma q+\delta p
\\
\gamma sx+\delta x
\end{pmatrix}.
$$
Passing to the variable $y=(\gamma s+\delta)x $
we get
$$
\begin{pmatrix} 
\alpha q+\beta p\\
(\alpha s+\beta)(\gamma s+\delta)^{-1} x
\end{pmatrix}
=
\begin{pmatrix}
a&b\\c&d
\end{pmatrix}
\begin{pmatrix}
\gamma q+\delta p
\\
 x
\end{pmatrix}
$$
This  implies the desired statement.
\hfill $\square$

\sm


{\bf \punct Non-algebraically closed fields.}
Now let $K$ be a non-algebraically closed infinite field .

\begin{proposition}
Let 
$$
w(\lambda)=\frac{u(\lambda)}{v(\lambda)}
$$
be a rational function
 on $\P K^1$ such that $v(0)\ne 0$. Then it  is a characteristic function of a certain 
 element of $\Coll_\infty(1)$.
\end{proposition}

{\sc Proof} by induction. Pass to the variable $s=\lambda^{-1}$.
 We say that degree of $w(s)$ is 
 $\deg v(s) $ (since $s=\infty$ is not a pole of $ w(s)$, we have $\deg u(s)\le \deg v(s)$).
  For functions of degree $1$ the statement is correct.
Assume that the statement is correct for functions of degree $<n$. Consider a function $w(s)$
of degree  
$n$. Take a linear fractional transformation
$$
\wt w(s):=
\frac{\alpha w(s)+\beta}{\gamma w(s) +\delta}
$$
such that $w(s)$ has a pole at some finite point $\sigma$ and a zero at some point $\tau$.
Then we can decompose $\wt w(s)$:
$$
\wt w(s)= \frac{s-\tau}{s-\sigma}\, y(s)
,$$
where
$y(s)$ is a rational function of degree $<n$. Both factors are characteristic functions, 
therefore $\wt w(s)$ also is a characteristic function.
\hfill $\square$


\section{Maps of Bruhat-Tits trees}

\COUNTERS

 Now  $\K$ is the $p$-adic field $\Q_p$ 
 and $\O\subset \K$ is the ring of integers.
 All considerations below can be automatically extended to
 arbitrary locally compact non-Archimedian fields
 (few  words must be changed).
 
 \sm

 {\bf\punct Colligations.}
Denote by $\frcoll_N(\alpha)$ the set of all matrices 
$\begin{pmatrix}a&b\\c&d \end{pmatrix}$ over $\K$ of size $\alpha+N$  defined
up to the equivalence
\begin{equation}
\begin{pmatrix}a&b\\c&d \end{pmatrix}\sim
\begin{pmatrix}1&0\\0&u \end{pmatrix}
\begin{pmatrix}a&b\\c&d \end{pmatrix}
\begin{pmatrix}1&0\\0&u \end{pmatrix}^{-1},
\qquad 
\text{where $u\in\GL(N,\O)$.}
\label{eq:conjugacy-class}
\end{equation}

We define $\circ$-product 
$$
\frcoll_{N_1}(\alpha)\times \frcoll_{N_2}(\alpha)
\to \frcoll_{N_1+N_2}(\alpha)
$$
by the same formula (\ref{eq:circ}).

 As above we define 
$\frcoll_\infty(\alpha)$ and the associative $\circ$-product on $\frcoll_\infty(\alpha)$.

\sm

{\bf \punct Bruhat-Tits buildings.}
Consider a linear space $\K^n$ over $\K$. A {\it lattice} $R$
in $\K^n$ is a compact $\O$-submodule in $\K^n$ such that
$\K\cdot R=\K^n$. In other words (see, e.g. \cite{Weil}, \cite{Ner-gauss}),  
in a certain basis $e_j\in\K^n$, a submodule $R$ has the form
$\oplus_j \O e_j$.  The space $\Lat_n$ of all lattices
is a homogeneous space,
$$
 \Lat_n\simeq
\GL(n,\K)/ \GL(n,\O).
$$

We intend to construct two simplicial complexes
$\BT_n$ and $\BT_n^*$.

\sm

1) Consider an oriented graph, whose vertices are lattices in $\K^n$. We draw arrow from a vertex 
$R$ to a vertex $T$ if $T\supset R\supset pT$.  If $k$ vertices pairwise are connected by arrows,
then we draw a simplex with such vertices. In this way we get a simplicial complex  $\BT_n$, all maximal simplices
have dimension $n$. The group $\GL(n,\K)$ acts transitively on the set of all maximal simplices
(and also on the set of simplices of each given dimension $j=0$, $1$, \dots, $n$).

\sm

2) Consider a non-oriented graph whose vertices are lattices defined by a dilatation,
$R\sim R'$ if $R=\lambda R'$ for some $\lambda\in\K^\times$. Denote
$$
\Lat_n^*:=\Lat_n/\K^\times.
$$
 We connect two vertices $R\not\sim T$ by an edge if for some $\lambda$ we have $pT\subset\lambda R\subset T$. If $k$ vertices pairwise are connected by edges,
then we draw a simplex with such vertices.
 We get a simplicial complex $\BT_n^*$,
dimensions of all maximal simplices are $n-1$. The projective linear group 
$$
\PGL(n,\K)=\GL(n,\K)/\K^\times
$$
 acts transitively 
on the set of all simplices of a given dimension $j=0$, $1$,\dots, $n-1$. 

\sm

We have a natural map 
$$
\BT_n(\K)\to \BT^*_n(\K)
,$$
we send a lattice (a vertex) to the corresponding equivalence class, this induces a map of graphs. Moreover, vertices of a $k$-dimensional simplex fall to  vertices
of a simplex of dimension $\le k$. 

These complexes are called '{\it Bruhat-Tits buildings}', see, e.g., \cite{Gar}, \cite{Ner-gauss}.
For $n=2$ the building $\BT_2(\K)$ is an infinite  tree, each vertex is an end of
$(p+1)$ edges.

\sm


{\bf \punct Construction of characteristic functions.}
Consider the space $\K^2=\K^1\oplus \K^1$.
For any lattice $R\subset \K^2$ consider the lattice
$$
R\otimes \O^N\subset \K^2\otimes \K^N=\K^N\oplus \K^N
.
$$
For a colligation
$g=\begin{pmatrix}a&b\\c&d \end{pmatrix}$
we write the equation
\begin{equation}
\begin{pmatrix}q\\y\end{pmatrix}=
\begin{pmatrix}a&b\\c&d \end{pmatrix}
\begin{pmatrix}p\\x\end{pmatrix}
.
\label{eq:char-3}
\end{equation}
Consider the set $\chi_g(R)$ of all $q\oplus p\in \K^\alpha\oplus \K^\alpha$ such that there are
$y\oplus x\in R\otimes \K^N$ satisfying the equation (\ref{eq:char-3}). 

\begin{proposition}
\label{pr:lattice}
{\rm a)} The sets $\chi_g(R)$ are lattices.

\sm

{\rm b)} If $R$, $T\in \Lat_2$ are connected by an arrow, then $\chi_g(R)$ and $\chi_g(T)$
are connected by an arrow or coincide.

\sm

{\rm c)} A lattice $\chi_g(R)$ depends on the conjugacy class containing $g$ and not on
$g$ itself.

\sm

{\rm d)} $\chi_g(\lambda R)=\lambda\, \chi_g(R)$ for $\lambda\in \K^\times$.
\end{proposition}

Proof is given in the next subsection.

Thus $\chi_g$ is a map 
$$\chi_g:\BT_2(\K)\to \BT_n(\K).$$
 Since it commutes with multiplications
by scalars, we get also a well-defined map
$$
\chi_g^*:\BT_2(\K)\to \BT_n^*(\K).
$$

\sm


{\bf \punct Reformulation of the definition.%
\label{ss:reformulation}}
Consider the space 
$$
H=\K^\alpha\oplus \K^\alpha \oplus \K^N\oplus \K^N
$$
consisting of vectors with coordinates $q$, $p$, $y$, $x$. Consider the following subspaces
and submodules
in $H$:

\sm

--- $G\subset H$ is the graph of $g$;

\sm

--- $U =0\oplus 0\oplus \K^N\oplus \K^N$;

\sm

--- $V=\K^\alpha\oplus \K^\alpha\oplus 0\oplus 0$;

\sm

--- the $\O$-submodule $S=\K^\alpha\oplus \K^\alpha\oplus (R\otimes \K^N)$.
 
 \sm
 
Consider the intersection
$G\cap S$ and its projection to $V$ along $U$. The result is $\chi_g(R)$.
 
\sm

{\sc Proof of Proposition \ref{pr:lattice}.}
The statements a), b), d) follow from new version of the definition, c) 
follows from $\GL(N,\O)$-invariance of $R\otimes \O^N$.
\hfill $\square$

\sm


{\bf\punct Products.} Now we wish to obtain an analog of Theorem \ref{th:prod-1}.
For this purpose, we need a definition of multiplication of lattices.

Let $S$, $T\subset \K^\alpha\oplus\K^\alpha$ be lattices. We define their {\it product}
$ST$ as  the set  of all $u\oplus w\in \K^\alpha\oplus \K^\alpha$ such that there is
$v\in\K^\alpha$ satisfying $u\oplus v\in S$, $v\oplus w\in T$. This is the usual product 
of relations (or multi-valued maps), see, e.g., \cite{Ner-book}.

\begin{theorem}
\label{th:product}
$$
\chi_{g\circ  h}(S)=\chi_{g}(S) \chi_h(S)
.
$$
\end{theorem}

{\sc Proof.} Let $g=\begin{pmatrix}a&b\\c&d\end{pmatrix}$, 
$h=\begin{pmatrix}\alpha&\beta\\\gamma&\delta\end{pmatrix}$.
Let $r=\chi_g(S)q$, $q=\chi_h(S) p$. Then there are $z$, $y$, $y'$, $x$ such that
$$
y\oplus x\in R\otimes \O^{N_1},\qquad y'\oplus x'\in R\otimes\O^{N_2}
$$
satisfying
$$
\begin{pmatrix}r\\y \end{pmatrix}=
\begin{pmatrix}a&b\\c&d\end{pmatrix}
\begin{pmatrix}q\\ x\end{pmatrix},\qquad
\begin{pmatrix}q\\y' \end{pmatrix}=
\begin{pmatrix}\alpha&\beta\\\gamma&\delta\end{pmatrix}
\begin{pmatrix}p\\x'\end{pmatrix}
.
$$
Then 
$$
\begin{pmatrix}r\\y\\y' \end{pmatrix}
=
\begin{pmatrix}a&b&0\\c&d&0\\0&0&1\end{pmatrix}
\begin{pmatrix}q\\ x\\y'\end{pmatrix}
=
\begin{pmatrix}a&b&0\\c&d&0\\0&0&1\end{pmatrix}
\begin{pmatrix}\alpha&0&\beta\\0&1&0\\ \gamma&0&\delta \end{pmatrix}
\begin{pmatrix}q\\ x\\x'\end{pmatrix}
.$$
This proves the desired statement.
\hfill $\square$

\sm

Consider  the natural projection 
$$
\mathrm{pr}:
\frcoll_N(\alpha)\to \Coll_n(\alpha).
$$
Formally, we have two characteristic functions of an element of $\frcoll_N(\alpha)$, one is defined 
on $\P \K^1$, another on $\BT_2(\K)$. In fact, the second function is the value of the first on the boundary
of the building. Now we intend to explain this.

\sm


{\bf\punct Convergence of lattices to subspaces.} We say that a sequence of lattices
$R_j\in \Lat_n$ converges  to a subspace $L\subset \K^n$ if 

\sm

a) For each $\epsilon$ for sufficiently large $j$ a lattice $R_j$ is contained in the $\epsilon$-neighborhood
of $L$

\sm

b) For each compact set $S\subset L$ we have $R_j\cap L\subset S$ for sufficiently large $j$.

\begin{proposition}
\label{pr:convergence}
Let a sequence $R^j\in\Lat_n$ converge to a subspace $L\subset\O^n$. Let
$M\subset \O^n$ be a subspace. Let $\pi:\O^n\to \O^n/L$ be the natural projection. 
Then 

\sm

{\rm a)} $R_j\cap M$ converges c $L\cap M$.

\sm

{\rm b)} $\pi(R_j)$ converges to $\pi (M)$.
\end{proposition}

The statement is obvious.

\sm

We say that a sequence $R_j^*\in \Lat^*$ converges to a subspace $L$ if we have
a convergence $R_j\to L$ for some representatives of $R_j^*$. Notice that a sequence 
$R_j^*$ can have many limits in this sense%
\footnote{Moreover, $0$ and $\O^n$ are limits of all sequences according our definition.}%
$^,$%
\footnote{See \cite{Ner-hausdorff}.}.
 However a limit subspace of a given dimension
is unique.

\sm


{\bf\punct Boundary values.}

\begin{proposition}
\label{pr:continuity}
Let $g\in \frcoll_N(\alpha)$.
Let $\lambda\in \P \K^1$ be a nonsingular point
of the characteristic function $\chi_{\mathrm{pr}(g)}(\lambda)$ defined on $\P\K^1$.  
Let $L$ be the line in $\K^2$ corresponding $\lambda$. Let $R_j\in\Lat_n(\K)$
converges to $\ell$. Then $\chi_g(R_j)$ converges to $\chi_{\mathrm{pr}(g)}(\lambda)$
\end{proposition}

{\sc Proof.} The statement follows from Subsection \ref{ss:reformulation}
and Proposition \ref{pr:convergence}.
\hfill $\square$

\sm


{\bf\punct Rational maps of Bruhat--Tits trees.}

\begin{corollary}
Any rational map $\P\K^1\to \P\K^1$ can be extended to a continuous map
of Bruhat-Tits trees, such that image of a vertex is a vertex and image of an
edge is an edge or a vertex.
\end{corollary}

{\sc Proof.}
Represent a rational map as a characteristic
function of a colligation $q\in\Coll_\infty(1)$. We take a colligation
$g\in\frcoll_\infty(1)$ such that $\mathrm{pr}(g)=q$, and take
the corresponding map $\BT_2^*(\K)\to \BT_2^*(\K)$.


\section{Rational maps of buildings.}

\COUNTERS

{\bf\punct $m$-colligations.}
Fix $\alpha\ge 0$, $m\ge 1$. Let $N>0$. Consider the space
$\Mat(\alpha+mN,\K)$
of block matrices of size $\alpha+N+\dots+N$. Denote by $\frcoll_N(\alpha|m)$
the set of such matrices up to the equivalence
\begin{multline}
\begin{pmatrix}
a&b_1&\dots&b_m\\
c_1&d_{11}&\dots&d_{1m}\\
\vdots&\vdots&\ddots&\vdots\\
c_m&d_{m1}&\dots&d_{mm}
\end{pmatrix}
\sim\\\sim
\begin{pmatrix}
1&0&\dots&0\\
0&u&\dots&0\\
\vdots&\vdots&\ddots&\vdots\\
0&0&\dots&u
\end{pmatrix}
\begin{pmatrix}
a&b_1&\dots&b_m\\
c_1&d_{11}&\dots&d_{1m}\\
\vdots&\vdots&\ddots&\vdots\\
c_m&d_{m1}&\dots&d_{mm}
\end{pmatrix}
\begin{pmatrix}
1&0&\dots&0\\
0&u&\dots&0\\
\vdots&\vdots&\ddots&\vdots\\
0&0&\dots&u
\end{pmatrix}^{-1},\\
\text{where $u\in \GL(N,\O)$}.
 \label{eq:equivalence-long}
\end{multline}

We define a multiplication
$$
\frcoll_{N_1}(\alpha|m)\times \frcoll_{N_2}(\alpha|m)
\to
\frcoll_{N_1+N_2}(\alpha|m)
$$
by
\begin{multline*}
\begin{pmatrix}
a&b_1&\dots&b_m\\
c_1&d_{11}&\dots&d_{1m}\\
\vdots&\vdots&\ddots&\vdots\\
c_m&d_{m1}&\dots&d_{mm}
\end{pmatrix}\circ
\begin{pmatrix}
p&q_1&\dots&q_m\\
r_1&t_{11}&\dots&t_{1m}\\
\vdots&\vdots&\ddots&\vdots\\
r_m&t_{m1}&\dots&t_{mm}
\end{pmatrix}
=\\=
\begin{pmatrix}
a&b_1&0&\dots&b_m&0\\
c_1&d_{11}&0&\dots&d_{1m}&0\\
0&0&1_{N_2}&\dots&0&0\\
\vdots&\vdots&\vdots&\ddots&\vdots&\vdots\\
c_m&d_{m1}&0&\dots&d_{mm}&0\\
0&0&0&\dots&0&1_{N_2}
\end{pmatrix}
\begin{pmatrix}
p&0&q_1&\dots&0&q_m\\
0&1_{N_1}&0&\dots&0&0\\
r_1&0&t_{11}&\dots&0&t_{1m}\\
\vdots&\vdots&\vdots&\ddots&\vdots&\vdots\\
0&0&0&\dots&1_{N_1}&0\\
r_m&t_{m1}&0&\dots&0&t_{mm}
\end{pmatrix}
=\\=
\begin{pmatrix}
a&b_1&aq_1&\dots&b_m&aq_m\\
c_1 p&d_{11}&c_1q_1&\dots&d_{1m}&c_1q_m\\
r_1&0&t_{11}&\dots&0&t_{1m}\\
\vdots&\vdots&\vdots&\ddots&\vdots&\vdots\\
c_m p&d_{m1}&c_mq_1&\dots&d_{mm}&c_mq_m\\
r_m&0&t_{m1}&\dots&0&t_{mm}\\
\end{pmatrix}
\end{multline*}


{\bf\punct Characteristic functions.} For a lattice
$R\in\Lat_{2m}(\K)$ consider the lattice
$$
R\otimes \O^N\subset \K^{2m}\otimes \K^{N}=
 (\K^{m}\otimes \K^{N}) \oplus (\K^{m}\otimes \K^{N})
$$
 For $g\in \Mat(\alpha+km)$ we write the following equation
 \begin{equation}
 \begin{pmatrix}
 q\\y_1\\\vdots\\y_m
 \end{pmatrix}=
 \begin{pmatrix}
a&b_1&\dots&b_m\\
c_1&d_{11}&\dots&d_{1m}\\
\vdots&\vdots&\ddots&\vdots\\
c_m&d_{m1}&\dots&d_{mm}
\end{pmatrix}
 \begin{pmatrix}
 p\\x_1\\\vdots\\x_m
 \end{pmatrix}
 ,
 \label{eq:char-last}
 \end{equation}
where $p$, $q$ range in $\K^\alpha$, and $x_j$, $y_j\in\K^N$.
Denote by $\chi_g(R)$ the set of all $q\oplus p\in\K^{2\alpha}$
such that there exists $y\oplus x\in R\otimes \O^N$,
for which  equality (\ref{eq:char-last}) holds.

\begin{theorem}
{\rm a)} $\chi_g(R)$ is a lattice in $\K^\alpha\oplus \K^\alpha$.

\sm

{\rm b)}  The characteristic function $\chi_g(R)$ is an invariant of the equivalence
{\rm(\ref{eq:equivalence-long})}.

\sm

{\rm c)} The map $\chi_g:\Lat_{2m}\to \Lat_{2\alpha}$ induces maps
$$
\BT_{2m}\to\BT_{2\alpha}, \qquad \BT^*_{2m}\to\BT^*_{2\alpha}.
$$

\sm

{\rm c)} For any $g\in \frcoll_{N_1}(\alpha|m)$,  $h\in \frcoll_{N_2}(\alpha|m)$,
the following identity holds
$$
\chi_{g\circ h}(R)=\chi_g(R)\chi_h(R).
$$
\end{theorem}

{\sc Proofs} repeats the proofs given above for $m=1$.
See, also a more sophisticated objects in \cite{Ner-buildings}.
\hfill $\square$

\sm


{\bf\punct Extension to the boundary.} Next (see \cite{Ner-char}, \cite{Ner-invariant}), we extend characteristic functions to
the distinguished  boundaries of buildings.
Let $S\in \Mat(m,\K)$. Again write equation (\ref{eq:char-last}).
We say 
 $q=\chi_g(S) p$ if there exists
$y$ such that $q$, $p$, $y$, $x=Sy$ satisfy the equation (\ref{eq:char-last}).
In other words,
$$
\chi_g(S)=a+b\wt S(1-d\wt S)^{-1} c
,$$
where $\wt S=S\otimes 1_N$,
$$
\wt S= \begin{pmatrix}
s_{11}\cdot 1_N&\dots &s_{1m}\cdot 1_N\\
\vdots&\ddots &\vdots\\
s_{m1}\cdot 1_N&\dots &s_{mm}\cdot 1_N
\end{pmatrix}
$$

\begin{theorem}
{\rm a)} For any $g\in\frcoll_{N_1}(\alpha|m)$, $h\in\frcoll_{N_2}(\alpha|m)$,
$$ 
\chi_{g\circ h}(S)=\chi_g(S)\chi_h(S).
$$

{\rm b)} If a sequence of lattices $R_j\in\Lat(\K^m\oplus\K^m)$ converges to the
graph of $S$, then $\chi_g(R_j)$ converges to $\chi_g(S)$.
\end{theorem}

Proof is  the same as above for $m=1$.

{\tt Math.Dept., University of Vienna,

 Nordbergstrasse, 15,
Vienna, Austria

\&

Institute for Theoretical and Experimental Physics,

Bolshaya Cheremushkinskaya, 25, Moscow 117259,
Russia

\&

Mech.Math.Dept., Moscow State University,

Vorob'evy Gory, Moscow

e-mail: neretin(at) mccme.ru

URL:www.mat.univie.ac.at/$\sim$neretin

wwwth.itep.ru/$\sim$neretin
}

\end{document}